\documentclass[12pt]{article}

\usepackage{latexsym}

\newtheorem{theorem}{Theorem}[section]

\newtheorem{example}[theorem]{Example}
\newtheorem{definition}[theorem]{Definition}

\newtheorem{proposition}[theorem]{Proposition}

\newtheorem{corollary}[theorem]{Corollary}

\newenvironment{proof}{\medskip\noindent{\it Proof.\ }}{\hfill \mbox{$\Box$}\medskip}

\begin{document}

\def\eqnsep{50pt}

\def\sof{\hfill\rule{2mm}{2mm}}
\def\ls{\leq}
\def\gs{\geq}
\def\SS{\mathcal S}
\def\qq{{\bold q}}
\def\txx{{\frac1{2\sqrt{x}}}}

\title{Restricted Permutations, Fibonacci Numbers, and $k$-generalized Fibonacci Numbers}
   
\author{
Eric S. Egge \\
Department of Mathematics \\
Gettysburg College\\
Gettysburg, PA  17325  USA \\
\\
eggee@member.ams.org \\
\\
Toufik Mansour \\
LaBRI (UMR 5800), Universit\'e Bordeaux 1, 351 cours de la Lib\'eration,\\
33405 Talence Cedex, France\\
\\
toufik@labri.fr
}

\maketitle

\begin{abstract}
A permutation $\pi \in S_n$ is said to {\it avoid} a permutation $\sigma \in S_k$ whenever $\pi$ contains no subsequence with all of the same pairwise comparisons as $\sigma$.
For any set $R$ of permutations, we write $S_n(R)$ to denote the set of permutations in $S_n$ which avoid every permutation in $R$.
In 1985 Simion and Schmidt showed that $|S_n(132, 213, 123)|$ is equal to the Fibonacci number $F_{n+1}$.
In this paper we generalize this result in several ways.
We first use a result of Mansour to show that for any permutation $\tau$ in a certain infinite family of permutations, $|S_n(132, 213, \tau)|$ is given in terms of Fibonacci numbers or $k$-generalized Fibonacci numbers.
In many cases we give explicit enumerations, which we prove bijectively.
We then use generating function techniques to show that for any permutation $\gamma$ in a second infinite family of permutations, $|S_n(123, 132, \gamma)|$ is also given in terms of Fibonacci numbers or $k$-generalized Fibonacci numbers.
In many cases we give explicit enumerations, some of which we prove bijectively.
We go on to use generating function techniques to show that for any permutation $\omega$ in a third infinite family of permutations, $|S_n(132, 2341, \omega)|$ is given in terms of Fibonacci numbers, and for any permutation $\mu$ in a fourth infinite family of permutations, $|S_n(132, 3241, \mu)|$ is given in terms of Fibonacci numbers and $k$-generalized Fibonacci numbers.
In several cases we give explicit enumerations.
We conclude by giving an infinite class of examples of a set $R$ of permutations for which $|S_n(R)|$ satisfies a linear homogeneous recurrence relation with constant coefficients.
\end{abstract}

\section{Introduction}

Let $S_n$ denote the set of permutations of $\{1, \ldots, n\}$, written in one-line notation, and suppose $\pi \in S_n$ and $\sigma \in S_k$.
We say $\pi$ {\it avoids} $\sigma$ whenever $\pi$ contains no subsequence with all of the same pairwise comparisons as $\sigma$.
For example, the permutation 214538769 avoids 312 and 2413, but it has 2586 as a subsequence so it does not avoid 1243.
If $\pi$ avoids $\sigma$ then $\sigma$ is sometimes called a {\it pattern} or a {\it forbidden subsequence} and $\pi$ is sometimes called a {\it restricted permutation} or a {\it pattern-avoiding permutation}. 
Pattern avoidance has proved to be a useful language in a variety of seemingly unrelated problems, from stack sorting \cite[Ch. 2.2.1]{KnuthVol1}, to singularities of Schubert varieties \cite{LakSong}, to Chebyshev polynomials of the second kind \cite{CW,Krattenthaler,MansourVainshtein2}, to rook polynomials for a rectangular board \cite{MansourVainshtein3}. 

One important and often difficult problem in the study of restricted permutations is the enumeration problem:  given a set $R$ of permutations, enumerate the set $S_n(R)$ consisting of those permutations in $S_n$ which avoid every element of $R$.
The earliest solution to an instance of this problem seems to be MacMahon's enumeration of $S_n(123)$, which is implicit in chapter V of \cite{MacMahon}.
The first explicit solution seems to be Hammersley's enumeration of $S_n(321)$ in \cite{Hammersley}.
In \cite[Ch. 2.2.1]{KnuthVol1} and \cite[Ch. 5.1.4]{KnuthVol3} Knuth shows that for any $\sigma \in S_3$, we have $|S_n(\sigma)| = C_n$, the $n$th Catalan number.
Other authors considered restricted permutations in the 1970s and early 1980s (see, for instance, \cite{Rogers}, \cite{Rotem1}, and \cite{Rotem2}) but the first systematic study was not undertaken until 1985, when Simion and Schmidt \cite{SimionSchmidt} solved the enumeration problem for every subset of $S_3$.
More recent work on various instances of the enumeration problem may be found in \cite{Atkinson}, \cite{BDPP}, \cite{BonaExact}, \cite{BonaSmooth}, \cite{G}, \cite{Kremer1}, \cite{KremerShiu}, \cite{Mansour33k}, \cite{Mproc}, \cite{MansourVainshtein}, \cite{MansourVainshtein3}, \cite{MansourVainshtein2}, \cite{Stankova1}, \cite{Stankova2}, \cite{StankovaWestHex}, \cite{WestCatalanSchroder}, and \cite{WestGenTrees}. 

In this paper we are concerned with instances of the enumeration problem whose solutions involve the Fibonacci numbers $F_n$ or the $k$-generalized Fibonacci numbers $F_{k,n}$.
(See section 2 for definitions.)
The earliest example of such a result is Simion and Schmidt's proof \cite[Prop. 15]{SimionSchmidt} that 
\begin{equation}
\label{eqn:SimionSchmidtintro}
|S_n(123, 132, 213)| = F_{n+1} \hspace{\eqnsep} (n \ge 1).
\end{equation}
Somewhat later West used generating trees to show \cite{WestGenTrees} that for many sets $R$ consisting of one pattern of length three and one of length four, $|S_n(R)| = F_{2n-1}$.
More recently Mansour expressed \cite[Thm. 3]{Mansour33k} the ordinary generating function for $|S_n(132, 213, \tau)|$ as a determinant for every $\tau \in S_k(132, 213)$.
Using this result, Mansour showed that
$$|S_n(132, 213, 2341)| = F_{n+2} - 1 \hspace{30pt} (n \ge 1)$$
and
$$|S_n(132, 213, 1234)| = T_{n+1} \hspace{30pt} (n \ge 1).$$
Here $T_n$ is the $n$th Tribonacci number, defined by $T_0 = 0$, $T_1 = T_2 = 1$, and $T_n = T_{n-1} + T_{n-2} + T_{n-3}$ for $n \ge 3$.
Mansour also expressed \cite[Thm. 1]{Mansour33k} the ordinary generating function for $|S_n(123, 132, \tau)|$ as a determinant for a large class of permutations $\tau \in S_k(123, 132)$.
Using this result, Mansour showed that
$$|S_n(123, 132, 3241)| = F_{n+2} - 1 \hspace{30pt} (n \ge 1)$$
and
$$|S_n(123, 132, 3214)| = T_{n+1} \hspace{30pt} (n \ge 1).$$

In this paper we begin by using Mansour's determinant formula \cite[Thm. 3]{Mansour33k} to show that if $a$, $b$, and $c$ are nonnegative integers with $b \ge 1$ and $a+c \ge 1$ then
$$\sum_{n=0}^\infty |S_n(132, 213, \beta_{a,b,c})| x^n = \frac{(1-x)^{a+c} +x^b \left( \sum\limits_{i=0}^{a+c-1} (1-x)^i x^{a+c-i+1} \right)}{(1-x)^{a+c} (1 - x - \ldots - x^{b-1})},$$
where $\beta_{a,b,c}$ is the permutation in $S_{a+b+c}$ given by
$$\beta_{a,b,c} = a+b+c, a+b+c-1, \ldots, b+c+1, c+1, c+2, \ldots, b+c, c, c-1, \ldots, 2,1.$$
We go on to give a bijective proof that for $n \ge 1$,
$$|S_n(132, 213, \beta_{a,b,c})| = \sum_{k=1}^{a+c-1} {{n-1} \choose {k-1}} + \sum_{k=a+c}^n {{k-1} \choose {a+c-1}} F_{b-1,n-k+1}.$$

We then turn our attention to some of those permutations $\gamma \in S_n(123, 132)$ to which Mansour's determinant formula \cite[Thm. 1]{Mansour33k} does not apply.
For all nonnegative integers $a$, $b$, and $c$, we give recurrence relations for the generating function for $|S_n(123, 132, \gamma_{a,b,c})|$, where $\gamma_{a,b,c}$ is the permutation in $S_{a+b+c+1}$ given by
$$\gamma_{a,b,c} = a+b+c+1, a+b+c, \ldots, b+c+2, b+c, b+c-1, \ldots, c+1, b+c+1, c, c-1, \ldots, 2,1.$$
Using these relations, we show that for all $n \ge 3$ and all $b \ge 3$,
\begin{equation}
\label{eqn:022intro}
|S_n(132, 123, \gamma_{0,2,2})| = F_{n+2} + F_n - 3
\end{equation}
and
\begin{equation}
\label{eqn:0b2intro}
|S_n(132, 123, \gamma_{0,b,2})| = \left( \frac{1}{b-1} \right) \left( 3 F_{b,n+1} + (b-4) F_{b,n-1} + 3 \sum_{i=3}^{b-1} (b-i) F_{b,n-i+1} - 3 \right).
\end{equation}
We give several other explicit enumerations;  for instance, we show that for all $n \ge 3$,
$$|S_n(132, 123, \gamma_{1,2,2})| = F_{n+3}+ F_{n+1} - 3n +2$$
and for all $n \ge 5$,
$$|S_n(132, 123, \gamma_{2,2,2})| = 5 F_{n+1} - 9n + 21.$$
We give bijective proofs of (\ref{eqn:022intro}) and (\ref{eqn:0b2intro}).
We also enumerate $S_n(132, 123, \gamma_{1,b,1})$ for all $b \ge 2$, obtaining a sum of $b$-generalized Fibonacci numbers plus a polynomial of degree one in $n$.
Finding a bijective proof of these enumerations is an open problem.

Next we prove two generalizations of Simion and Schmidt's result (\ref{eqn:SimionSchmidtintro}) which are not related to Mansour's determinant formulas.
For all $k \ge 4$, we give a recurrence relation for the generating function for $|S_n(132, 2341, \omega_k)|$, where $\omega_k$ is the permutation in $S_k$ given by 
$$\omega_k = k, k-1, \ldots, 4, 2, 1, 3.$$
Using this relation we give several explicit enumerations.
For example, we show that for all $n \ge 2$,
$$|S_n(132, 2341, \omega_4)| = F_{n+5} - {{n+1} \choose {2}} - 2 {{n+1} \choose {1}} - 2.$$
For all nonnegative integers $a$ and $b$, we give recurrence relations for the generating function for $|S_n(132, 3241, \mu_{a,b})|$, where $\mu_{a,b}$ is the permutation in $S_{a+b}$ given by
$$\mu_{a,b} = b+a, b+a-1, \ldots, b+1, 1, 2, \ldots, b.$$
Using these relations we give several explicit enumerations.
For example, we show that for all $n \ge 1$,
$$|S_n(132, 3241, \mu_{2,3})| = F_{n+8} - {{n+1}\choose{4}} - 3 {{n+1} \choose {3}} - 4 {{n+1} \choose {2}} - 9 {{n+1} \choose {1}} - 11.$$
It is an open problem to find bijective proofs of these enumerations.

We conclude our main results by expanding our point of view to include all sequences which satisfy a linear homogeneous recurrence relation with constant coefficients, instead of restricting our attention to those sequences whose terms can be expressed in terms of Fibonacci numbers or $k$-generalized Fibonacci numbers.
More specifically, let $k$ and $l$ denote positive integers such that $l \le k$ and let $a_1, \ldots, a_l$ denote a sequence of positive integers in which $a_1, \ldots, a_{l-1}$ are distinct.
For a certain set $R_{a_1, \ldots, a_l}^k \subseteq S_k$, we show that for all $n \ge k$,
$$|S_n(R^k_{a_1, \ldots, a_k})| = \sum_{j=1}^l (k - a_j - \eta_j) |S_{n-j}(R^k_{a_1, \ldots, a_l})|,$$
where $\eta_j = \{i\ |\ a_i > a_j\ \mbox{and}\ i < j\}$ for $1 \le j \le l$.
(See section 7 for a precise definition of $R_{a_1, \ldots, a_l}^k$.)
This result includes Simion and Schmidt's result (\ref{eqn:SimionSchmidtintro}) as a special case, and we give several other examples of specific enumerations which arise from this result.
Among these examples is an infinite family of enumerations involving the Fibonacci numbers.

\section{Background and Notation}

Let $S_n$ denote the set of permutations of $[n] = \{1, \ldots, n\}$, written in one-line notation, and suppose $\pi \in S_n$ and $\sigma \in S_k$.
We say a permutation $\pi$ {\it avoids} a permutation $\sigma$ whenever $\pi$ contains no subsequence with all of the same pairwise comparisons as $\sigma$.
For example, the permutation 214538769 avoids 312 and 2413, but it has 2586 as a subsequence so it does not avoid 1243.
We make this idea precise in the following definition.

\begin{definition}
For any permutation $\pi \in S_n$ and any $i$ $(1 \le i \le n)$, we write $\pi(i)$ to denote the element of $\pi$ in position $i$.
We say a permutation $\pi \in S_n$ {\upshape avoids} a permutation $\sigma \in S_k$ whenever there is no sequence $1 \le i_{\sigma(1)} < i_{\sigma(2)} < \cdots < i_{\sigma(k)} \le n$ such that $\pi(i_1) < \pi(i_2) < \cdots < \pi(i_k)$.
\end{definition}

\noindent
If $\pi$ avoids $\sigma$ then $\pi$ is sometimes called a {\it restricted permutation} or a {\it pattern-avoiding permutation} and $\sigma$ is sometimes called a {\it forbidden subsequence}.
In this paper we will be interested in permutations which avoid several patterns, so for any set $R$ of permutations we write $S_n(R)$ to denote the elements of $S_n$ which avoid every element of $R$.
For any set $R$ of permutations we take $S_n(R)$ to be the empty set whenever $n < 0$ and we take $S_0(R)$ to be the set containing only the empty permutation.
When $R = \{\pi_1, \pi_2, \ldots, \pi_r\}$ we often write $S_n(R) = S_n(\pi_1, \pi_2, \dots, \pi_r)$.

For all integers $k \ge 1$, the $k$-generalized Fibonacci number $F_{k,n}$ satisfies the recurrence obtained by adding more terms to the recurrence for the Fibonacci numbers.
More specifically, we set $F_{k,n} = 0$ for all $n \le 0$ and we set $F_{k,1} = 1$.
For all $n \ge 2$ we define $F_{k,n}$ recursively by setting
\begin{equation}
\label{eqn:kgenrecurrence}
F_{k,n} = \sum_{i=1}^k F_{k,n-i} \hspace{\eqnsep} (n \ge 2).
\end{equation}
The term ``$k$-generalized Fibonacci number'' alludes to the fact that $F_{2,n} = F_n$ for all $n$.
It is routine to verify that the ordinary generating function for $F_{k,n}$ is given by
\begin{equation}
\label{eqn:Fkngf}
\sum_{n=0}^\infty F_{k,n} x^n = \frac{x}{1 - x - x^2 - \ldots - x^k} \hspace{30pt} (k \ge 1).
\end{equation}
We will also make use of the following combinatorial interpretation of $F_{k,n}$.

\begin{proposition}
\label{prop:kgencombinatorial}
The number of tilings of a $1 \times n$ rectangle with tiles of size $1 \times 1$, $1 \times 2$, \ldots, $1 \times k$ is the $k$-generalized Fibonacci number $F_{k,n+1}$.
\end{proposition}
\begin{proof}
The result is immediate for $n \le 1$, so it suffices to show that the number of such tilings satisfies (\ref{eqn:kgenrecurrence}).
To do this, observe there is a one-to-one correspondence between tilings of a $1 \times (n - i)$ rectangle and tilings of a $1 \times n$ rectangle in which the rightmost tile has length $i$.
Therefore, if we count tilings of a $1 \times n$ rectangle according to the length of the rightmost tile, we find the number of such tilings satisfies (\ref{eqn:kgenrecurrence}), as desired.
\end{proof}

For more information concerning the $k$-generalized Fibonacci numbers, see \cite{Flores}, \cite{Gabai}, \cite{Lynch}, \cite{Miles}, and \cite{Miller}.

\section{Generalizations of Simion and Schmidt Involving 132 and 213}
\label{sec:MansourDetCor}

In this section we generalize Simion and Schmidt's result (\ref{eqn:SimionSchmidtintro}) by replacing 123 with a longer permutation.
One natural way to do this is to replace 123 with $123\ldots k$ for some $k \ge 4$.
We go further, however, replacing 123 with the permutation given in the following definition.

\begin{definition}
\label{defn:tausubr}
For all positive integers $k$, $m$, and $r_0, \ldots, r_m$ such that $k+1 = r_0 > r_1 > \cdots > r_m = 1$, we write $\tau_{r_0, \ldots, r_m}$ to denote the permutation in $S_k$ given by
$$\tau_{r_0, \ldots, r_m} = r_1, r_1 + 1, \ldots, k, r_2, r_2 + 1, \ldots, r_1 - 1, \ldots, 1, 2, \ldots r_{m-1} - 1.$$
\end{definition}

In this section we consider $|S_n(132, 213, \tau_{r_0, \ldots, r_m})|$ for certain sequences $r_0 > r_1 > \cdots > r_m$.
Since $\tau_{4,1} = 123$, these results generalize Simion and Schmidt's result (\ref{eqn:SimionSchmidtintro}).
We begin with a result of Mansour concerning the ordinary generating function for $|S_n(132, 213, \tau_{r_0, \ldots, r_m})|$.

\begin{theorem}
\label{thm:MansourDet}
(\cite[Thm. 3]{Mansour33k})
Fix positive integers $k, m$, and $r_0, \ldots, r_m$ such that $k+1 = r_0 > r_1 > \cdots > r_m = 1$ and let $\tau_{r_0, \ldots, r_m}$ be as in Definition \ref{defn:tausubr}.
Then
\begin{equation}
\label{eqn:MansourDet}
\sum_{n=0}^\infty |S_n(132, 213, \tau_{r_0, \ldots, r_m})| x^n = \det\left( 
\matrix{f_{r_0 - r_1} & -g_{r_0 - r_1} & 0 & \cdots & 0 \cr
f_{r_1 - r_2} & 1 & -g_{r_1 - r_2} & \ddots & 0 \cr
f_{r_2 - r_3} & 0 & 1 & \ddots & 0 \cr
\vdots & \vdots & \ddots & \ddots & -g_{r_{m-2} - r_{m-1}} \cr
f_{r_{m-1} - r_m} & 0 & 0 & \cdots & 1}
\right)
\end{equation}
where
${\displaystyle f_i(x) = \frac{1-x}{1- 2x + x^i}}$ and ${\displaystyle g_i(x) = \frac{x^i}{1 - 2x + x^i}}$ for $i \ge 0$.
\end{theorem}

In \cite{Mansour33k} Mansour focuses on the case of this theorem in which $k = 4$.
Here we use the result to give explicit enumerations involving Fibonacci numbers and $k$-generalized Fibonacci numbers.
Many such enumerations can be obtained from this theorem;  we consider just a few of them.
We begin by considering the case in which 123 is replaced with its most natural lengthening, $12\ldots k$.

\begin{corollary}
For all integers $n$ and all $k \ge 2$,
\begin{equation}
\label{eqn:12kenumeration}
|S_n(12\ldots k, 132, 213)| = F_{k-1,n+1}.
\end{equation}
Moreover,
\begin{equation}
\label{eqn:12kgf}
\sum_{n=0}^\infty |S_n(12\ldots k, 132, 213)| x^n = \frac{1}{1 - x - x^2 - \cdots - x^{k-1}}.
\end{equation}
\end{corollary}
\begin{proof}
Since $12\ldots k = \tau_{k+1,1}$, we set $m = 1$ and $r_1 = 1$ in (\ref{eqn:MansourDet}) to obtain (\ref{eqn:12kgf}).
Line (\ref{eqn:12kenumeration}) is immediate from (\ref{eqn:12kgf}), in view of (\ref{eqn:Fkngf}).
\end{proof}

Throughout this section we will give bijective proofs of our explicit enumerations.
All of these proofs are based on a certain constructive bijection between $S_n(132, 213)$ and the set of tilings of a $1 \times n$ rectangle with tiles of size $1 \times i$, where $i \ge 1$.
To describe this bijection, we first make an observation concerning the form of a permutation in $S_n(132, 213)$.

\begin{proposition}
\label{prop:2nform}
Fix $n \ge 0$ and suppose $\pi \in S_n$.
Then $\pi \in S_n(132, 213)$ if and only if there exists a positive integer $m$ and a sequence $n+1 = r_0 > r_1 > \cdots > r_m = 1$ such that $\pi = \tau_{r_0, r_1, \ldots, r_m}$.
\end{proposition}
\begin{proof}
It is routine to verify that $\tau_{r_0, r_1, \ldots, r_m}$ avoids 132 and 213.
To prove the converse, suppose $\pi \in S_n(132, 213)$.
Since $\pi$ avoids 132, every element to the left of $n$ in $\pi$ is larger than every element to the right of $n$.
Since $\pi$ avoids 213, the elements to the left of $n$ are in increasing order, so there exists $r_1$ such that $\pi = r_1, r_1+1, \ldots, n, \tilde{\pi}$, where $\tilde{\pi} \in S_{r_1 - 1}(132, 213)$.
By induction on $n$ there exists a positive integer $m$ and a sequence $n+1 = r_0 > r_1 > \cdots > r_m = 1$ such that $\pi = \tau_{r_0, \ldots, r_m}$, as desired.
\end{proof}

We now describe our bijection between $S_n(132, 213)$ and the set of tilings of a $1 \times n$ rectangle with tiles of size $1 \times i$, where $i \ge 1$.

\begin{definition}
\label{defn:calF}
For all $n \ge 0$, we write ${\cal F}_n$ to denote the map from $S_n(132, 213)$ to the set of tilings of a $1 \times n$ rectangle with tiles of size $1 \times i$, where $i \ge 1$, for which ${\cal F}_n(\tau_{r_0, \ldots, r_m})$ is the tiling whose tiles, when read from left to right, have sizes $1 \times (r_0-r_1), 1 \times (r_1 - r_2), \ldots, 1 \times (r_{m-1} - r_m)$.
\end{definition}

\begin{example}
We have
$${\cal F}_9(978652341) = \Box\ \Box\hspace{-2.5pt}\Box\ \Box\ \Box\ \Box\hspace{-2.5pt}\Box\hspace{-2.5pt}\Box\ \Box$$
and 
$${\cal F}_8(56782341) = \Box\hspace{-2.5pt}\Box\hspace{-2.5pt}\Box\hspace{-2.5pt}\Box\hspace{-2.5pt}\ \Box\hspace{-2.5pt}\Box\hspace{-2.5pt}\Box\ \Box.$$
\end{example}

\begin{proposition}
\label{prop:Fnbijection}
For all $n \ge 0$, the map ${\cal F}_n$ is a constructive bijection between $S_n(132, 213)$ to the set of tilings of a $1 \times n$ rectangle with tiles of size $1 \times i$, where $i \ge 1$.
\end{proposition}
\begin{proof}
To show ${\cal F}_n$ is a bijection, we construct ${\cal F}_n^{-1}$.
To do this, let $G_n$ denote the map from the set of tilings of a $1 \times n$ rectangle with tiles of size $1 \times i$, where $i \ge 1$, to $S_n(132, 213)$, which is given as follows.
For a given tiling, fill the rightmost tile with the numbers $1,2,\ldots$ from left to right, up to the length of the tile.
Fill the tile immediately to the left of the rightmost tile from left to right in the same way, beginning with the smallest available number.
Repeat this process until every tile is filled.
It is routine to verify that ${\cal F}_n(G_n(\mu)) = \mu$ for any tiling $\mu$, and that $G_n({\cal F}_n(\pi)) = \pi$ for any $\pi \in S_n(132,213)$.
Therefore $G_n = {\cal F}_n^{-1}$, so ${\cal F}_n$ is a bijection, as desired.
\end{proof}

The function ${\cal F}_n$ gives us a bijective proof of (\ref{eqn:12kenumeration}).

\begin{theorem}
\label{thm:12kbijection}
For all $n \ge 1$ and all $k \ge 2$, the restriction of ${\cal F}_n$ to the set \\ $S_n(12\ldots k, 132, 213)$ is a bijection between $S_n(12\ldots k, 132, 213)$ and the set of tilings of a $1 \times n$ rectangle with tiles of size $1 \times 1$, $1 \times 2$, \ldots, $1 \times (k-1)$.
\end{theorem}
\begin{proof}
Observe that $\tau_{r_0, \ldots, r_m}$ avoids $12\ldots k$ if and only if $r_{i-1} - r_i \le k-1$ for $1 \le i \le m$.
Now the result is immediate from Propositions \ref{prop:2nform} and \ref{prop:Fnbijection}.
\end{proof}

Next we consider a more complicated lengthening of 123, in which we increase every entry in $12\ldots k$ by a fixed amount and then append a decreasing sequence to the result.
This longer permutation is $\tau_{s+t+2, s+1, s, \ldots, 2,1}$ for positive integers $s$ and $t$, but for notational convenience we abbreviate it as follows.

\begin{definition}
For all positive integers $s$ and $t$, we write $\alpha_{s,t}$ to denote the permutation in $S_{s+t+1}$ given by
\begin{displaymath}
\alpha_{s,t} = \tau_{s+t+2, s+1, s, s-1, \ldots, 2,1}.
\end{displaymath}
\end{definition}

\begin{proposition}
\label{prop:sm}
Let $s$ and $t$ denote positive integers.
Then
\begin{equation}
\label{eqn:smgf}
\sum_{n=0}^\infty |S_n(132, 213, \alpha_{s,t})| x^n = \frac{(1-x)^s + x^{t+1}\left( \sum\limits_{i=0}^{s-1} (1-x)^i x^{s-i-1} \right)}{(1-x)^s (1 - x - \cdots - x^t)}.
\end{equation}
\end{proposition}
\begin{proof}
In (\ref{eqn:MansourDet}), set $m=s+2$, $r_0 = s+t+2$, and $r_i = s+2-i$ for $1 \le i \le m-1$ to find
$$\sum_{n=0}^\infty |S_n(132, 213, \alpha_{s,t})| x^n = \det\left(
\matrix{f_{t+1}(x) & -g_{t+1}(x) & 0 & \cdots & 0 & 0 \cr
1 & 1 & - \frac{x}{1-x} & \cdots & 0 & 0 \cr
1 & 0 & 1 & \cdots & 0 & 0 \cr
\vdots & \vdots & \vdots & \ddots & \vdots & \vdots \cr
1 & 0 & 0 & \cdots & 1 & -\frac{x}{1-x} \cr
1 & 0 & 0 & \cdots & 0 & 1 \cr
}\right),$$
where the matrix on the right is $s+1$ by $s+1$.
To evaluate this determinant, let
$$A_n(x) = \det\left(
\matrix{f_{t+1}(x) & -g_{t+1}(x) & 0 & \cdots & 0 & 0 \cr
1 & 1 & - \frac{x}{1-x} & \cdots & 0 & 0 \cr
1 & 0 & 1 & \cdots & 0 & 0 \cr
\vdots & \vdots & \vdots & \ddots & \vdots & \vdots \cr
1 & 0 & 0 & \cdots & 1 & -\frac{x}{1-x} \cr
1 & 0 & 0 & \cdots & 0 & 1 \cr
}\right) \hspace{30pt} (n \ge 2),$$
where the matrix on the right is $n \times n$.
Expand this determinant along the bottom row to find
$$A_n(x) = A_{n-1}(x) + g_{t+1}(x) \left(\frac{x}{1-x}\right)^{n-2} \hspace{30pt} (n \ge 3).$$
Iterate this relation and use the fact that $A_2(x) = f_{t+1}(x) + g_{t+1}(x)$ to obtain
$$A_n(x) = f_{t+1}(x) + g_{t+1}(x) \sum_{i=0}^{n-2} \left( \frac{x}{1-x} \right)^i \hspace{30pt} (n \ge 2).$$
Since ${\displaystyle f_{t+1}(x) = \frac{1}{1 - x -\cdots - x^t}}$ and ${\displaystyle g_{t+1}(x) = \frac{x^{t+1}}{(1-x)(1-x-\cdots-x^t)}}$, we have
$$A_n(x) = \frac{(1-x)^{n-1} + x^{t+1} \sum\limits_{i=0}^{n-2} x^i (1-x)^{n-2-i}}{(1-x)^{n-1}(1-x-\cdots-x^t)} \hspace{30pt} (n \ge 2).$$
Set $n = s+1$ in the last line to obtain (\ref{eqn:smgf}), as desired.
\end{proof}

Rather than extract an explicit enumeration of $S_n(132, 213, \alpha_{s,t})$ directly from (\ref{eqn:smgf}), we use ${\cal F}_n$ to obtain one bijectively.

\begin{theorem}
\label{thm:smbijection}
Let $s$ and $t$ denote positive integers.
Then the restriction of ${\cal F}_n$ to the set $S_n(132, 213, \alpha_{s,t})$ is a constructive bijection between $S_n(132, 213, \alpha_{s,t})$ and the set of tilings of a $1 \times n$ rectangle in which all tiles except for the rightmost $s$ tiles have length at most $t$. 
In particular, for all $n \ge 1$,
\begin{equation}
\label{eqn:sm}
|S_n(132, 213, \alpha_{s,t})| = \sum_{k=1}^{s-1} {{n-1}\choose{k-1}} + \sum_{k=s}^n {{k-1}\choose{s-1}} F_{t,n-k+1}.
\end{equation} 
\end{theorem}
\begin{proof}
In view of Propositions \ref{prop:2nform} and \ref{prop:Fnbijection}, to prove the first part of the theorem it is sufficient to show that $\tau_{r_0, \ldots, r_m}$ avoids $\alpha_{s,t}$ if and only if $r_{i-1} - r_i < t+1$ for $1 \le i \le m-s$.
It is routine to verify that if $r_{i-1} - r_i < t+1$ for $1 \le i \le m-s$ then $\tau_{r_0, \ldots, r_m}$ avoids $\alpha_{s,t}$.
To show the converse, observe that if there exists $i$ with $1 \le i \le m-s$ and $r_{i-1} - r_i \ge t+1$ then the subsequence $r_i, r_i+1, \ldots, r_{i-1}-1, r_{i+1}, r_{i+2}, \ldots, r_{i+s}$ forms a pattern of type $\alpha_{s,t}$.

Line (\ref{eqn:sm}) follows from the first part of the theorem by a routine counting argument.
\end{proof}

We conclude this section by using (\ref{eqn:MansourDet}) and ${\cal F}_n$ to obtain the following general result, which is analogous to Proposition \ref{prop:sm} and Theorem \ref{thm:smbijection}.
For notational convenience, in this result we abbreviate
$$\beta_{a,b,c} = \tau_{a+b+c+1, a+b+c, \ldots, b+c+1, c+1, c, \ldots, 2,1}$$
for all nonnegative integers $a,b,$ and $c$ such that $a+c \ge 1$ and $b \ge 1$.

\begin{theorem}
Let $a$ and $c$ denote nonnegative integers such that $a+c \ge 1$ and let $b$ denote a positive integer.
Then the restriction of ${\cal F}_n$ to the set $S_n(132, 213, \beta_{a,b,c})$ is a constructive bijection between $S_n(132, 213, \beta_{a,b,c})$ and the set of tilings of a $1 \times n$ rectangle in which all tiles except for the leftmost $a$ tiles and the rightmost $c$ tiles have length at most $b-1$.
In particular, for $n \ge 1$,
\begin{equation}
\label{eqn:tauabc}
|S_n(132, 213, \beta_{a,b,c})| = \sum_{k=1}^{a+c-1} {{n-1} \choose {k-1}} + \sum_{k=a+c}^n {{k-1} \choose {a+c-1}} F_{b-1,n-k+1}.
\end{equation}
Moreover, 
\begin{equation}
\label{eqn:tauabcgf}
\sum_{n=0}^\infty |S_n(132, 213, \beta_{a,b,c})| x^n = \frac{(1-x)^{a+c} +x^b \left( \sum\limits_{i=0}^{a+c-1} (1-x)^i x^{a+c-i+1} \right)}{(1-x)^{a+c} (1 - x - \ldots - x^{b-1})}.
\end{equation}
\end{theorem}
\begin{proof}
In view of Propositions \ref{prop:2nform} and \ref{prop:Fnbijection}, to prove the first part of the theorem it is sufficient to show that $\tau_{r_0, \ldots, r_m}$ avoids $\beta_{a,b,c}$ if and only if $r_{i-1} - r_i < b$ for $a+1 \le i \le m-c$.
It is routine to verify that if $r_{i-1} - r_i < b$ for $a+1 \le i \le m-c$ then $\tau_{r_0, \ldots, r_m}$ avoids $\beta_{a,b,c}$.
To show the converse, observe that if there exists $i$ with $a+1 \le i \le m-c$ and $r_{i-1} - r_i \ge b$ then the subsequence
$$r_1, r_2, \ldots, r_a, r_i, r_i+1, \ldots, r_i+b-1, r_{m-c}+1, r_{m-c}+2, \ldots, r_m$$
forms a pattern of type $\beta_{a,b,c}$.
It follows that $\tau_{r_0, \ldots, r_m}$ avoids $\beta_{a,b,c}$ if and only if $r_{i-1} - r_i < b$ for $a+1 \le i \le m-c$, so the first part of the theorem holds.

Line (\ref{eqn:tauabc}) follows from the first part of the theorem by a routine counting argument.

To obtain (\ref{eqn:tauabcgf}), observe that if we set $s = a+c$ and $t = b-1$ in the expression on the right side of (\ref{eqn:sm}) then we obtain the expression on the right side of (\ref{eqn:tauabc}).
Therefore, in view of (\ref{eqn:smgf}), line (\ref{eqn:tauabcgf}) follows from (\ref{eqn:tauabc}).
\end{proof}

\section{Generalizations of Simion and Schmidt Involving 123 and 132}
\label{sec:MoreSS}

In section \ref{sec:MansourDetCor} we generalized Simion and Schmidt's result (\ref{eqn:SimionSchmidtintro}) by replacing 123 with a longer permutation.
In this section we generalize (\ref{eqn:SimionSchmidtintro}) by replacing 213 with a longer permutation.
The generalizations we obtain in this way are equivalent to those we would obtain by replacing 132 with a longer permutation, since 132 is the image of 213 under the reverse-complement map and 123 is fixed by this map.
We begin by setting some notation.

\begin{definition}
\label{defn:fabc}
For all nonnegative integers $a$, $b$, and $c$, we write $\gamma_{a,b,c}$ to denote the permutation in $S_{a+b+c+1}$ given by
\begin{displaymath}
\gamma_{a,b,c} = a+b+c+1, a+b+c, \ldots, b+c+2, b+c, b+c-1, \ldots, c+1, b+c+1, c, c-1, \ldots, 2,1
\end{displaymath}
We write $f_{a,b,c}(x)$ to denote the generating function
\begin{displaymath}
f_{a,b,c}(x) = \sum_{n=0}^\infty |S_n(123, 132, \gamma_{a,b,c})| x^n.
\end{displaymath}
\end{definition}

In this section we give a recursive formula for $f_{a,b,c}(x)$, which we use to give explicit enumerations of $S_n(123, 132, \gamma_{a,b,c})$ for certain values of $a$, $b$, and $c$.
Since $\gamma_{0,2,0} = 213$, these results generalize Simion and Schmidt's result (\ref{eqn:SimionSchmidtintro}).
We begin by considering the case in which $a=b=0$.

\begin{proposition}
Let $c$ denote a positive integer.
Then $f_{0,0,0}(x) = 1$ and 
\begin{equation}
\label{eqn:00c}
f_{0,0,c}(x) = 1 + x f_{0,0,c-1}(x) + \sum_{r=2}^{c+1} x^r f_{0,0,c-r+1}(x).
\end{equation}
\end{proposition}
\begin{proof}
To see that $f_{0,0,0}(x) = 1$, first observe that $\gamma_{0,0,0} = 1$.
Since the empty permutation is the only permutation (of any length) which avoids $1$, we have $f_{0,0,0}(x) = 1$.

To prove (\ref{eqn:00c}), first observe that the empty permutation avoids 123, 132, and $\gamma_{0,0,c}$.
For $n \ge 1$, observe that every permutation $\pi \in S_n(123, 132, \gamma_{0,0,c})$ for which $\pi(1) = n$ has the form $n, \tilde{\pi}$, where $\tilde{\pi} \in S_n(123, 132, \gamma_{0,0,c})$, and the resulting map is a bijection between $S_n(123, 132, \gamma_{0,0,c})$ and the set of permutations in $S_n(123, 132, \gamma_{0,0,c})$ which begin with $n$.
Similarly, observe that if $2 \le i < n$ then every permutation $\pi \in S_n(123, 132, \gamma_{0,0,c})$ for which $\pi(i) = n$ has the form $n-1, n-2, \ldots, n-i+1, n, \tilde{\pi}$, where $\tilde{\pi} \in S_{n-i}(123, 132, \gamma_{0,0,c-i+1})$, and the resulting map is a bijection between $S_{n-i}(123, 132, \gamma_{0,0,c-i+1})$ and the set of permutations in $S_n(123, 132, \gamma_{0,0,c})$ in which $n$ is in the $i$th position.
Combine these observations to obtain (\ref{eqn:00c}).
\end{proof}

We now consider $f_{0,b,0}(x)$.

\begin{proposition}
For any positive integer $b$,
\begin{equation}
\label{eqn:f0b0gf}
f_{0,b,0}(x) = \frac{1}{1-x-\cdots-x^b}.
\end{equation}
In particular,
\begin{equation}
\label{eqn:f0b0}
|S_n(123, 132, \gamma_{0,b,0})| = F_{b,n+1} \hspace{30pt} (n \ge 1).
\end{equation}
\end{proposition}
\begin{proof}
Line (\ref{eqn:f0b0gf}) is immediate from \cite[Thm. 1]{Mansour33k}, but for completeness we include a proof here.

First observe that the empty permutation avoids 123, 132, and $\gamma_{0,b,0}$.
Now observe that for $1 \le i \le b$, every permutation $\pi \in S_n(123, 132, \gamma_{0,b,0})$ with $\pi(i) = n$ has the form $n-1, n-2, \ldots, n-i+1, n, \tilde{\pi}$, where $\tilde{\pi} \in S_{n-i}(123, 132, \gamma_{0,b,0})$, and the resulting map is a bijection between $S_{n-i}(123, 132, \gamma_{0,b,0})$ and the set of permutations in $S_n(123, 132, \gamma_{0,b,0})$ in which $n$ is in the $i$th position.
Finally, observe that if $\pi$ avoids 123 and 132 and has $\pi(i) = n$ for $b+1 \le i \le n$ then $\pi$ does not avoid $\gamma_{0,b,0}$.
Combine these observations to obtain
$$f_{0,b,0}(x) = 1 + f_{0,b,0}(x) \sum_{r=1}^b x^r.$$
Solve this equation for $f_{0,b,0}(x)$ to obtain (\ref{eqn:f0b0gf}).

Line (\ref{eqn:f0b0}) is immediate from (\ref{eqn:f0b0gf}), in view of (\ref{eqn:Fkngf}).
\end{proof}

We also give a bijective proof of (\ref{eqn:f0b0}).

\begin{theorem}
\label{thm:kk1bijection}
For all $n \ge 1$ and all $b \ge 2$ there exists a constructive bijection between $S_n(123, 132, \gamma_{0,b,0})$ and the set of tilings of a $1 \times n$ rectangle with tiles of size $1 \times 1$, $1 \times 2$, \ldots, $1 \times b$.
\end{theorem}
\begin{proof}
Suppose we are given such a tiling.
We construct its corresponding permutation as follows.
Let $m$ denote the length of the rightmost tile.
Fill this tile with the numbers $1, 2, \ldots, m$ from left to right in the order $m-1, m-2, \ldots, 1, m$.
Fill the tile immediately to the left of the rightmost tile in the same way, using the smallest available numbers.
Repeat this process until every tile is filled.
It is routine to verify that this map is invertible, and that the permutation constructed avoids 123, 132, and $\gamma_{0,b,0}$.
\end{proof}

\begin{example}
Under the bijection given in the proof of Theorem \ref{thm:kk1bijection}, the permutation \\ 876954231 corresponds to the tiling $\Box\hspace{-2.5pt}\Box\hspace{-2.5pt}\Box\hspace{-2.5pt}\Box\ \Box\ \Box\ \Box\hspace{-2.5pt}\Box\ \Box$ and the permutation 986743512 corresponds to the tiling $\Box\ \Box\ \Box\hspace{-2.5pt}\Box\ \Box\hspace{-2.5pt}\Box\hspace{-2.5pt}\Box\ \Box\hspace{-2.5pt}\Box$.
\end{example}

We now find $f_{0,b,c}(x)$ in terms of $f_{0,0,c}(x)$.

\begin{proposition}
For any positive integers $b$ and $c$, 
\begin{equation}
\label{eqn:fbc}
f_{0,b,c}(x) = \frac{1 - x + x^{b+1} f_{0,0,c-1}(x)}{(1-x)(1-x-\cdots - x^b)}.
\end{equation}
\end{proposition}
\begin{proof}
First observe that the empty permutation avoids 123, 132, and $\gamma_{0,b,c}$.
Now observe that for $1 \le i \le b$, every permutation $\pi \in S_n(123, 132, \gamma_{0,b,c})$ with $\pi(i) = n$ has the form $n-1, n-2, \ldots, n-i+1, n, \tilde{\pi}$, where $\tilde{\pi} \in S_{n-i}(123, 132, \gamma_{0,b,c})$, and the resulting map is a bijection between $S_{n-i}(123, 132, \gamma_{0,b,c})$ and the set of permutations in $S_n(123, 132, \gamma_{0,b,c})$ in which $n$ is in the $i$th position.
Similarly, observe that for $b+1 \le i \le n$, every permutation $\pi \in S_n(123, 132, \gamma_{0,b,c})$ with $\pi(i) = n$ has the form $n-1, n-2, \ldots, n-i+1, n, \tilde{\pi}$, where $\tilde{\pi} \in S_{n-i}(123, 132, \gamma_{0,0,c-1})$, and the resulting map is a bijection between $S_{n-i}(123, 132, \gamma_{0,0,c-1})$ and the set of permutations in $S_n(123, 132, \gamma_{0,b,c})$ in which $n$ is in the $i$th position.
Combine these observations to find
$$f_{0,b,c}(x) = 1 + \sum_{i=1}^b x^i f_{0,b,c}(x) + \frac{x^{b+1}}{1-x} f_{0,0,c-1}(x).$$
Solve this equation for $f_{0,b,c}(x)$ to obtain (\ref{eqn:fbc}).
\end{proof}

We now obtain our recurrence relation for $f_{a,b,c}(x)$.

\begin{proposition}
For any positive integers $a$ and $b$ and any nonnegative integer $c$,
\begin{equation}
\label{eqn:fabc}
f_{a,b,c}(x) = 1 + x f_{a-1, b, c}(x) + \sum_{r=2}^a x^r f_{a-r+1,b,c}(x) + \frac{x^{a+1}}{1-x} f_{0,b,c}(x).
\end{equation}
\end{proposition}
\begin{proof}
First observe that the empty permutation avoids 123, 132, and $\gamma_{a,b,c}$.
Now observe that every permutation $\pi \in S_n(123, 132, \gamma_{a,b,c})$ for which $\pi(1) = n$ has the form $n, \tilde{\pi}$, where $\tilde{\pi} \in S_{n-1}(123, 132, \gamma_{a-1,b,c})$, and the resulting map is a bijection between $S_{n-1}(123, 132, \gamma_{a-1,b,c})$ and the set permutations in $S_n(123, 132, \gamma_{a,b,c})$ which begin with $n$.
Similarly, observe that for $2 \le i \le a$, every permutation $\pi \in S_n(123, 132, \gamma_{a,b,c})$ with $\pi(i) = n$ has the form $n-1, n-2, \ldots, n-i+1, n, \tilde{\pi}$, where $\tilde{\pi} \in S_{n-i}(123, 132, \gamma_{a-i+1,b,c})$, and the resulting map is a bijection between $S_{n-i}(123, 132, \gamma_{a-i+1,b,c})$ and the set of permutations in $S_n(123, 132, \gamma_{a,b,c})$ in which $n$ is in the $i$th position.
Finally, observe that for $a+1 \le i \le n$, every permutation $\pi \in S_n(123, 132, \gamma_{a,b,c})$ for which $\pi(i) = n$ has the form $n-1, n-2, \ldots, n-i+1, n, \tilde{\pi}$, where $\tilde{\pi} \in S_{n-i}(123, 132, \gamma_{0,b,c})$, and the resulting map is a bijection between $S_{n-i}(123, 132, \gamma_{0,b,c})$ and the set of permutations in $S_n(123, 132, \gamma_{a,b,c})$ in which $n$ is in the $i$th position.
Combine these observations to obtain (\ref{eqn:fabc}).
\end{proof}

We now use our recurrence relations (\ref{eqn:00c}), (\ref{eqn:fbc}), and (\ref{eqn:fabc}) to obtain explicit enumerations of $S_n(132, 123, \gamma_{a,b,c})$ for various values of $a$, $b$, and $c$.
We begin with the case in which $a = 0$, $c = 1$, and $b$ is arbitrary.

\begin{proposition}
Let $b$ denote a positive integer.
Then 
\begin{equation}
\label{eqn:0b1}
|S_n(132, 123, \gamma_{0,b,1})| = \sum_{k=1}^n F_{b,k} \hspace{30pt} (n \ge 1).
\end{equation}
Moreover,
\begin{equation}
\label{eqn:0b1gf}
\sum_{n=0}^\infty |S_n(132, 123, \gamma_{0,b,1})| x^n = 1 + \frac{x}{1-2x+x^{b+1}}.
\end{equation}
\end{proposition}
\begin{proof}
To prove (\ref{eqn:0b1gf}), set $c = 1$ in (\ref{eqn:fbc}), obtaining
$$\sum_{n=0}^\infty |S_n(132, 123, \gamma_{0,b,1})| x^n = \frac{1-x+x^{b+1}}{(1-x)(1-x-\cdots-x^b)} = 1 + \frac{x}{1-2x+x^{b+1}},$$
as desired.
To prove (\ref{eqn:0b1}) observe that in view of (\ref{eqn:Fkngf}),
\begin{eqnarray*}
\frac{x}{1-2x+x^{b+1}} &=& \left( \sum_{n=0}^\infty x^n \right) \left( \sum_{n=0}^\infty F_{b,n} x^n \right) \\
&=& \sum_{n=0}^\infty \left( \sum_{k=1}^n F_{b,k} \right) x^n.
\end{eqnarray*}
Now (\ref{eqn:0b1}) is immediate from (\ref{eqn:0b1gf}).
\end{proof}

\noindent
{\bf Remark}
Line (\ref{eqn:0b1gf}) is also immediate from \cite[Thm. 1]{Mansour33k}.

\medskip

We also give a bijective proof of (\ref{eqn:0b1}).

\begin{theorem}
\label{thm:0b1bijection}
Let $b$ denote a positive integer.
There exists a constructive bijection between $S_n(132, 123, \gamma_{0,b,1})$ and the set of tilings of a $1 \times (n+1)$ rectangle with tiles of size $1 \times i$ $(i \ge 1)$ in which all tiles except the rightmost tile have length at most $b$, and the rightmost tile has length at least 2.
\end{theorem}
\begin{proof}
Suppose we are given such a tiling.
We construct its corresponding permutation as follows.
Let $m$ denote the length of the rightmost tile.
Fill this tile with the numbers $1, 2, \ldots, m$ from left to right in the order $m-2, m-3, \ldots, 1, m-1$, leaving the rightmost square empty.
Fill the tile immediately to the left of the rightmost tile using the smallest available numbers and in the same pattern as the rightmost tile, but put a number in every square.
Repeat this process until every square of every remaining tile is filled.
It is routine to verify that this map is invertible, and that the permutation constructed avoids 123, 132, and $\gamma_{0,b,1}$.
\end{proof}

\begin{example}
Fix $b = 3$.
Under the bijection given in the proof of Theorem \ref{thm:0b1bijection}, the permutation 879653214 corresponds to the tiling $\Box\hspace{-2.5pt}\Box\hspace{-2.5pt}\Box\ \Box\ \Box\ \Box\hspace{-2.5pt}\Box\hspace{-2.5pt}\Box\hspace{-2.5pt}\Box\hspace{-2.5pt}\Box$ and the permutation $976845231$ corresponds to the tiling $\Box\ \Box\hspace{-2.5pt}\Box\hspace{-2.5pt}\Box\ \Box\hspace{-2.5pt}\Box\ \Box\hspace{-2.5pt}\Box\ \Box\hspace{-2.5pt}\Box$.
\end{example}

Next we enumerate $S_n(132, 123, \gamma_{a,b,c})$ when $a = 0$, $c = 2$, and $b$ is arbitrary.
We first consider the case in which $b = 2$.

\begin{proposition}
For all $n \ge 3$,
\begin{equation}
\label{eqn:022}
|S_n(132, 123, \gamma_{0,2,2})| = F_{n+2} + F_n - 3.
\end{equation}
Moreover,
\begin{equation}
\label{eqn:022gf}
\sum_{n=0}^\infty |S_n(132, 123, \gamma_{0,2,2})| x^n = 3 + x + x^2 + \frac{x^2+4x-2}{1-2x+x^3}.
\end{equation}
\end{proposition}
\begin{proof}
To obtain (\ref{eqn:022gf}), set $c = 2$ in (\ref{eqn:fbc}) and use (\ref{eqn:00c}) to obtain
\begin{eqnarray*}
\sum_{n=0}^\infty |S_n(132, 123, \gamma_{0,2,2})| x^n &=& \frac{1-x+x^3+x^4+x^5}{(1-x)(1-x-x^2)} \\
&=& 3 + x + x^2 + \frac{x^2 + 4x - 2}{1-2x+x^3},
\end{eqnarray*}
as desired.
To obtain (\ref{eqn:022}), observe that
$$\frac{x^2 + 4x - 2}{1-2x+x^3} = \frac{-3}{1-x} + \frac{2x + 1}{1-x-x^2}.$$
Now (\ref{eqn:022}) is immediate in view of (\ref{eqn:Fkngf}).
\end{proof}

We now handle the case in which $b \ge 3$.

\begin{proposition}
Let $b$ denote a positive integer such that $b \ge 3$.
Then for all $n \ge 3$,
\begin{equation}
\label{eqn:0b2}
|S_n(132, 123, \gamma_{0,b,2})| = \left( \frac{1}{b-1} \right) \left( 3 F_{b,n+1} + (b-4) F_{b,n-1} + 3 \sum_{i=3}^{b-1} (b-i) F_{b,n-i+1} - 3 \right).
\end{equation}
Alternatively, for all $n \ge 3$,
\begin{equation}
\label{eqn:0b2alt}
|S_n(132, 123, \gamma_{0,b,2})| = \sum_{k=1}^{n-1} F_{b,k} + 2 \sum_{k=1}^{n-2} F_{b,k}.
\end{equation}
Moreover, 
\begin{equation}
\label{eqn:0b2gf}
\sum_{n=0}^\infty |S_n(132, 123, \gamma_{0,b,2})| x^n = 1 + x + x^2 + \frac{x^2 + 2x^3}{1-2x+x^{b+1}}.
\end{equation}
\end{proposition}
\begin{proof}
To obtain (\ref{eqn:0b2gf}), set $c = 2$ in (\ref{eqn:fbc}) and use (\ref{eqn:00c}) to obtain
\begin{eqnarray*}
\sum_{n=0}^\infty |S_n(132, 123, \gamma_{0,b,2})| x^n &=& \frac{1 - x + x^{b+1} + x^{b+2}+x^{b+3}}{(1-x)(1-x-\ldots-x^b)} \\
&=& 1 + x + x^2 + \frac{x^2 + 2x^3}{1-2x+x^{b+1}},
\end{eqnarray*}
as desired.
To obtain (\ref{eqn:0b2}), observe that if $b \ge 3$ then
$$\frac{x^2 + 2x^3}{1-2x+x^{b+1}} = \left( \frac{1}{b-1} \right) \left( \frac{-3}{1-x} + \frac{3 + (b-4) x^2 + 3 \sum\limits_{i=3}^{b-1} (b-i) x^i}{1-x-\ldots -x^b} \right).$$
Now (\ref{eqn:0b2}) is immediate in view of (\ref{eqn:Fkngf}).
To obtain (\ref{eqn:0b2alt}), observe that
\begin{eqnarray*}
\frac{x^2 + 2 x^3}{(1-x)(1-x-\ldots-x^b)} &=& \left( \sum_{n=0}^\infty x^n \right) \left( \sum_{n=0}^\infty \left( F_{b,n-1} + 2 F_{b,n-2} \right) x^n \right) \\
&=& \sum_{n=0}^\infty \left( \sum_{k=0}^{n-1} F_{b,k} + 2 \sum_{k=0}^{n-2} F_{b,k} \right) x^n. \\
\end{eqnarray*}
Now (\ref{eqn:0b2alt}) is immediate in view of (\ref{eqn:Fkngf}).
\end{proof}

We also give a bijective proof of (\ref{eqn:022}) and (\ref{eqn:0b2alt}).

\begin{theorem}
\label{thm:0220b2bijection}
Let $b \ge 2$.
For all $n \ge 1$, let $B_n$ denote the set of tilings of a $1 \times n$ rectangle for which the following hold.
\begin{enumerate}
\item There is at most one tile with length greater than $b$;  call this tile the {\upshape long tile}.
\item There is at most one tile to the right of the long tile.
\item The tile to the right of the long tile has length 1 or 2.
\end{enumerate}
Then there exists a constructive bijection between $B_n$ and $S_n(132, 132, \gamma_{0,b,2})$.
\end{theorem}
\begin{proof}
Suppose we are given such a tiling.
We construct its corresponding permutation as follows.
Let $m$ denote the length of the rightmost tile.
Fill this tile with the numbers $1,2, \ldots, m$ from left to right in the order $m-1, m-2, \ldots, 1, m$.
Fill the tile immediately to the left of the rightmost tile in the same way, using the smallest available numbers.
Repeat this process until every tile is filled.
It is routine to verify that this map is invertible, and that the permutation constructed avoids 132, 123, and $\gamma_{0,b,2}$.
\end{proof}

\begin{example}
Fix $b = 3$.
Under the bijection given in Theorem \ref{thm:0220b2bijection}, the permutation  879643251 corresponds to the tiling $\Box\hspace{-2.5pt}\Box\hspace{-2.5pt}\Box\ \Box\ \Box\hspace{-2.5pt}\Box\hspace{-2.5pt}\Box\hspace{-2.5pt}\Box\ \Box$, the permutation $978643215$ corresponds to the tiling $\Box\ \Box\hspace{-2.5pt}\Box\ \Box\ \Box\hspace{-2.5pt}\Box\hspace{-2.5pt}\Box\hspace{-2.5pt}\Box\hspace{-2.5pt}\Box$, and the permutation $6543712$ corresponds to the tiling $\Box\hspace{-2.5pt}\Box\hspace{-2.5pt}\Box\hspace{-2.5pt}\Box\hspace{-2.5pt}\Box\ \Box\hspace{-2.5pt}\Box$.
\end{example}

Next we enumerate $S_n(132, 123, \gamma_{a,b,c})$ when $a = 1$, $c = 1$, and $b$ is arbitrary.

\begin{proposition}
Let $b$ denote a positive integer such that $b \ge 2$.
Then for all $n \ge 1$,
\begin{equation}
\label{eqn:1b1}
|S_n(132, 123, \gamma_{1,b,1})| = \left( \frac{1}{b-1} \right)^2 \left( \frac{b^2 - 3b}{2} + (1-b)n + g_b(n)\right),
\end{equation}
where
$$g_b(n) = \frac{b^2 - b + 2}{2} F_{b,n+1} + (b-1) F_{b,n} + \sum_{i=2}^{b-1} \left[ (1-b) {{i} \choose {2}} + \frac{b^2-b-2}{2} i - \frac{b^2-3b}{2} \right] F_{b,n-i+1}.$$
Moreover,
\begin{equation}
\label{eqn:1b1gf}
\sum_{n=0}^\infty |S_n(132, 123, \gamma_{1,b,1})| x^n = \frac{1 - 2x + x^2 + x^{b+1}}{(1-x)^2 (1- x - \ldots - x^b)}.
\end{equation}
\end{proposition}
\begin{proof}
To obtain (\ref{eqn:1b1gf}), set $a = c = 1$ in (\ref{eqn:fabc}) and use (\ref{eqn:0b1gf}) to obtain
$$\sum_{n=0}^\infty |S_n(132, 123, \gamma_{1,b,1})| x^n = \frac{1}{1-x} + \frac{x^2}{1-2x+x^{b+1}} + \frac{x^3}{1 - 3x + 2x^2 + x^{b+1} - x^{b+2}}.$$
Now (\ref{eqn:1b1gf}) is immediate.
To obtain (\ref{eqn:1b1}), first observe that
$$\frac{x^2}{1-2x+x^{b+1}} = \left( \frac{1}{b-1} \right) \left( \frac{-1}{1-x} + \frac{1 + \sum\limits_{i=2}^{b-1} (b-i) x^i}{1 - x - \ldots - x^b}\right)$$
and
$$\frac{x^3}{1 - 3x + 2x^2 + x^{b+1} - x^{b+2}} = \left( \frac{1}{b-1} \right)^2 \left( \frac{-{{b-2}\choose{2}}}{1-x} + \frac{1-b}{(1-x)^2} + \frac{h_b(x)}{1 - x -\ldots - x^b}\right),$$
where
$$h_b(x) = \frac{b^2 -3b + 4}{2} + (b-1) x + \sum\limits_{i=2}^{b-1} \left[ (i-2)(b-1)^2 - \left( {{i} \choose {2}} - 1 \right) (b-1) - (i-1) {{b-2} \choose {2}} \right] x^i.$$
Now (\ref{eqn:1b1}) is immediate in view of (\ref{eqn:Fkngf}).
\end{proof}

We conclude this section with several other explicit enumerations of $S_n(132, 123, \gamma_{a,b,c})$ for small values of $a$, $b$, and $c$.

\begin{proposition}
\label{prop:122}
For all $n \ge 3$,
\begin{equation}
\label{eqn:122}
|S_n(132, 123, \gamma_{1,2,2})| = F_{n+3}+ F_{n+1} - 3n +2.
\end{equation}
Moreover,
\begin{equation}
\label{eqn:122gf}
\sum_{n=0}^\infty |S_n(132, 123, \gamma_{1,2,2})| x^n = -4 -2x - x^2 + \frac{6x^3+4x^2-12x+5}{(1-x)^2 (1-x-x^2)}.
\end{equation}
\end{proposition}
\begin{proof}
To obtain (\ref{eqn:122gf}), set $a = 1$, $b = 2$, and $c = 2$ in (\ref{eqn:fabc}) and use (\ref{eqn:022gf}).
To obtain (\ref{eqn:122gf}), observe that
$$\frac{6x^3+4x^2-12x+5}{(1-x)^2 (1-x-x^2)} = \frac{5}{1-x} - \frac{3}{(1-x)^2} + \frac{x+3}{1-x-x^2}.$$
Now (\ref{eqn:122gf}) is immediate in view of (\ref{eqn:Fkngf}).
\end{proof}

We omit the proofs of the remaining propositions in this section, since they are similar to the proof of Proposition \ref{prop:122}.

\begin{proposition}
For all $n \ge 5$,
\begin{equation}
\label{eqn:222}
|S_n(132, 123, \gamma_{2,2,2})| = 5 F_{n+1} - 9n + 21.
\end{equation}
Moreover,
\begin{displaymath}
\sum_{n=0}^\infty |S_n(132, 123, \gamma_{2,2,2})| x^n = -25 - 16x - 11x^2 - 5 x^3 - 2x^4 + \frac{30x^3+14x^2-61x+26}{(1-x)^2(1-x-x^2)}.
\end{displaymath}
\end{proposition}

\begin{proposition}
For all $n \ge 7$,
\begin{equation}
\label{eqn:322}
|S_n(132, 123, \gamma_{3,2,2})| = 3 F_{n+2} + 3 F_n - 24n +91.
\end{equation}
Moreover,
\begin{eqnarray*}
\lefteqn{\sum_{n=0}^\infty |S_n(132, 123, \gamma_{3,2,2})| x^n = } \\
 & &  -93 - 75x - 53x^2-36x^3-20x^4-9x^5-2x^6 + \frac{121 x^3 + 15 x^2 - 206 x + 94}{(1-x)^2(1-x-x^2)}.
\end{eqnarray*}
\end{proposition}

\begin{proposition}
For all $n \ge 3$,
\begin{equation}
|S_n(132, 123, \gamma_{1,3,2})| = \frac{1}{2} \left( F_{3,n+2} + 2 F_{3,n} + F_{3,n-1} - 3n + 5 \right).
\end{equation}
Moreover,
\begin{displaymath}
\sum_{n=0}^\infty |S_n(132, 123, \gamma_{1,3,2})| x^n = -2 - 2x - x^2 + \frac{5x^4 + x^3 - 6x + 3}{(1-x)^2 (1 - x - x^2 - x^3)}.
\end{displaymath}
\end{proposition}

\begin{proposition}
For all $n \ge 5$,
\begin{equation}
|S_n(132, 123, \gamma_{2,3,2})| = \frac{1}{2} \left( 4 F_{3,n+1} + F_{3,n} - 3 F_{3,n-1} - 9n + 30 \right).
\end{equation}
Moreover,
\begin{eqnarray*}
\lefteqn{\sum_{n=0}^\infty |S_n(132, 123, \gamma_{2,3,2})| x^n = } \\
& & -16 -12x - 7x^2 - 5x^3-2x^4 + \frac{18x^4 + 8x^3+4x^2-38x+17}{(1-x)^2 (1-x-x^2-x^3)}.
\end{eqnarray*}
\end{proposition}

\begin{proposition}
For all $n \ge 7$,
\begin{equation}
|S_n(132, 123, \gamma_{3,3,2})| = \frac{1}{2}\left( 3 F_{3,n} + 10 F_{3,n-1} - 3 F_{3,n-2} - 24n + 115 \right).
\end{equation}
Moreover,
\begin{eqnarray*}
\lefteqn{\sum_{n=0}^\infty |S_n(132, 123, \gamma_{3,3,2})| x^n = } \\
& & -55 -46x - 38x^2 - 24x^3 - 16x^4-9x^5-2x^6 + \frac{76x^4 + 2x^3 + 11x^2 - 121x + 56}{(1-x)^2(1-x-x^2-x^3)}.
\end{eqnarray*}
\end{proposition}

\section{Generalizations of Simion and Schmidt Involving Two Longer Permutations}

In sections \ref{sec:MansourDetCor} and \ref{sec:MoreSS} we generalized Simion and Schmidt's result (\ref{eqn:SimionSchmidtintro}) by replacing one of the three forbidden subsequences with a longer forbidden subsequence.
In this section we generalize (\ref{eqn:SimionSchmidtintro}) by replacing two of the forbidden subsequences with longer forbidden subsequences in several ways.
We begin by setting some notation.

\begin{definition}
\label{defn:gk}
For all $k \ge 4$, we write $\omega_k$ to denote the permutation in $S_k$ given by 
$$\omega_k = k, k-1, \ldots, 5, 4, 2, 1, 3,$$
and we write $g_k(x)$ to denote the generating function
\begin{displaymath}
g_k(x) = \sum_{n=0}^\infty |S_n(132, 2341, \omega_k)| x^n.
\end{displaymath}
\end{definition}

As in the previous section, we give a recurrence relation for $g_k(x)$.

\begin{proposition}
We have 
\begin{equation}
\label{eqn:g3}
g_3(x) = \frac{x^3-x+1}{(1-x)(1-x-x^2)}.
\end{equation}
Moreover, for all $k > 3$,
\begin{equation}
\label{eqn:gk}
g_k(x) = \left( \frac{1}{1-x} \right) \left( 1 + x ( g_{k-1}(x) - 1 ) + \sum_{r=2}^{k-2} x^r ( g_{k-r+1}(x) - 1 ) + \frac{x^{k-1}}{1-x} ( g_3(x) - 1) \right).
\end{equation}
\end{proposition}
\begin{proof}
Using the notation of section 3, observe that $\omega_3 = 213$ and $\tau_{5,2,1} = 2341$.
Now (\ref{eqn:g3}) is immediate from (\ref{eqn:MansourDet}).

To obtain (\ref{eqn:gk}), first observe that the empty permutation avoids 132, 2341, and $\omega_k$.
For $n \ge 1$, observe that every permutation $\pi \in S_n(132, 2341, \omega_k)$ for which $\pi(n) = n$ has the form $\tilde{\pi}, n$, where $\tilde{\pi} \in S_{n-1}(132, 2341, \omega_k)$, and the resulting map is a bijection between $S_{n-1}(132, 2341, \omega_k)$ and the set of permutations in $S_n(132, 2341, \omega_k)$ which end with $n$.
Similarly, observe that if $1 \le i \le k-2$ and $i < n$ then every permutation $\pi \in S_n(132, 2341, \omega_k)$ with $\pi(i) = n$ has the form $n-1, n-2, \ldots, n-i+1, n, \tilde{\pi}$, where $\tilde{\pi} \in S_{n-i}(132, 2341, \omega_{k-i+1})$, and the resulting map is a bijection between $S_{n-i}(132, 2341, \omega_{k-i+1})$ and the set of permutations in $S_n(132, 2341, \omega_k)$ in which $n$ is in the $i$th position.
Finally, observe that if $k-1 \le i \le n-1$ then every permutation $\pi \in S_n(132, 2341, \omega_k)$ with $\pi(i) = n$ has the form $n-1, n-2, \ldots, n-i+1, n, \tilde{\pi}$, where $\tilde{\pi} \in S_{n-i}(132, 2341, \omega_3)$, and the resulting map is a bijection between $S_{n-i}(132, 2341, \omega_3)$ and the set of permutations in $S_n(132, 2341, \omega_k)$ in which $n$ is in the $i$th position.
Combine these observations to find
$$g_k(x) = 1 + x g_k(x) + x (g_{k-1}(x)-1) + \sum_{r=2}^{k-2} x^r (g_{k-r+1}(x) - 1) + \frac{x^{k-1}}{1-x} \left( g_3(x) - 1 \right).$$
Solve this equation for $g_k(x)$ to obtain (\ref{eqn:gk}).
\end{proof}

We now use (\ref{eqn:g3}) and (\ref{eqn:gk}) to obtain explicit enumerations of $S_n(132, 2341, \omega_k)$ for $k \le 6$.
We will use the following identities, whose proofs we omit.
\begin{equation}
\label{eqn:fnpluskgf}
\sum_{n=0}^\infty F_{n+k} x^n = \frac{F_{k-1} x + F_k}{1-x-x^2} \hspace{30pt} (k \ge 2)
\end{equation}
\begin{equation}
\label{eqn:nplus1choosekgf}
\sum_{n=0}^\infty {{n+1} \choose {k}} x^n = \frac{x^{max(k-1,0)}}{(1-x)^{k+1}} \hspace{30pt} (k \ge 0)
\end{equation}
We begin with $S_n(132, 2341, \omega_3)$.

\begin{proposition}
For all $n \ge 1$,
\begin{equation}
\label{eqn:omega3}
|S_n(132, 2341, \omega_3)| = F_{n+2} - 1.
\end{equation}
\end{proposition}
\begin{proof}
This result is proved in \cite{Mansour33k}, but for completeness we give a proof here.
Observe that
$$g_3(x) = 1 - \frac{1}{1-x} + \frac{F_1 x + F_2}{1 - x - x^2}.$$
Now (\ref{eqn:omega3}) is immediate from (\ref{eqn:fnpluskgf}) and (\ref{eqn:nplus1choosekgf}).
\end{proof}

We also give a bijective proof of (\ref{eqn:omega3}).

\begin{theorem}
\label{thm:Mansourfn+2-1}
For all $n \ge 1$, there exists a constructive bijection between $S_n(132, 213, 2341)$ and the set of tilings of a $1 \times (n+1)$ rectangle with tiles of size $1 \times 1$ and $1 \times 2$ using at least one $1 \times 2$ tile.
\end{theorem}
\begin{proof}
Suppose we are given such a tiling.
We construct its corresponding permutation as follows.
Replace the rightmost $1 \times 2$ tile with a 1, and fill the (necessarily $1 \times 1$) tiles to the right of the 1 with $2, 3, \ldots$ from left to right.
Now fill the rightmost empty tile with the smallest numbers available, placing them in the tile from left to right in increasing order.
Repeat this process until every tile is filled.
It is routine to verify that this map is invertible, and that the permutation constructed avoids 132, 213, and 2341, as desired.
\end{proof}

\begin{example}
Under the bijection given in the proof of Theorem \ref{thm:Mansourfn+2-1}, the permutation \\ 87564123 corresponds to the tiling $\Box\ \Box\ \Box\hspace{-2.5pt}\Box\ \Box\ \Box\hspace{-2.5pt}\Box\ \Box\ \Box$ and the permutation 86745321 corresponds to the tiling $\Box\ \Box\hspace{-2.5pt}\Box\ \Box\hspace{-2.5pt}\Box\ \Box\ \Box\ \Box\hspace{-2.5pt}\Box$.
\end{example}

\begin{proposition}
\label{prop:omega4}
For all $n \ge 1$,
\begin{equation}
\label{eqn:omega4}
|S_n(132, 2341, \omega_4)| = F_{n+5} - {{n+1} \choose {2}} - 2 {{n+1} \choose {1}} - 2.
\end{equation}
Moreover,
\begin{equation}
\label{eqn:omega4gf}
\sum_{n=0}^\infty |S_n(132, 2341, \omega_4)| = \frac{1 - 3x+3x^2+x^3-x^4}{(1 - x - x^2) (1 - x)^3}.
\end{equation}
\end{proposition}
\begin{proof}
To obtain (\ref{eqn:omega4gf}), set $k = 4$ in (\ref{eqn:gk}) and use (\ref{eqn:g3}) to simplify the result.
To obtain (\ref{eqn:omega4}), observe that
$$\frac{x^4 - x^3 - 3x^2 + 3x - 1}{(x^2 + x - 1) (1 - x)^3} = \frac{F_4 x + F_5}{1 - x - x^2} - \frac{2}{1-x} - \frac{2}{(1-x)^2} - \frac{x}{(1-x)^3}.$$
Now (\ref{eqn:omega4}) is immediate from (\ref{eqn:fnpluskgf}) and (\ref{eqn:nplus1choosekgf}).
\end{proof}

We omit the proofs of the next two propositions, since they are similar to the proof of Proposition \ref{prop:omega4}.

\begin{proposition}
For all $n \ge 2$,
\begin{equation}
\label{eqn:omega5}
|S_n(132, 2341, \omega_5)| = 3 F_{n+5} - 2 {{n+1} \choose {3}} - 4 {{n+1} \choose {2}} - 3 {{n+1} \choose {1}} - 14.
\end{equation}
Moreover,
\begin{displaymath}
\sum_{n=0}^\infty |S_n(132, 2341, \omega_5)| x^n = \frac{1 - 4x + 6x^2-2x^3-x^4+3x^5-x^7}{(1 - x - x^2)(1-x)^4}.
\end{displaymath}
\end{proposition}

\begin{proposition}
For all $n \ge 2$,
\begin{equation}
\label{eqn:omega6}
|S_n(132, 2341, \omega_6)| = 5 F_{n+6} + F_{n+4} - 4 {{n+1} \choose {4}} - 7 {{n+1} \choose {3}} - 5 {{n+1} \choose {2}} - 27 {{n+1} \choose {1}} - 8.
\end{equation}
Moreover,
\begin{displaymath}
\sum_{n=0}^\infty |S_n(132, 2341, \omega_k)| x^n = \frac{1 - 5x+10x^2-8x^3+x^4+5x^5+x^6+x^7-2x^8}{(1 - x - x^2)(1 - x)^5}.
\end{displaymath}
\end{proposition}

We now consider another way of replacing two of the forbidden subsequences in Simion and Schmidt's result (\ref{eqn:SimionSchmidtintro}) with longer permutations.
We begin by setting some notation.

\begin{definition}
For all nonnegative integers $a$ and $b$ we write $\mu_{a,b}$ to denote the permutation in $S_{a+b}$ given by
$$\mu_{a,b} = b+a, b+a-1, \ldots, b+1, 1, 2, \ldots, b$$
and we write $h_{a,b}(x)$ to denote the generating function
$$h_{a,b}(x) = \sum_{n=0}^\infty |S_n(132, 3241, \mu_{a,b})| x^n.$$
\end{definition}

We now give a recurrence relation for $h_{0,b}(x)$.

\begin{proposition}
Let $b$ denote a positive integer.
Then
\begin{equation}
\label{eqn:h01}
h_{0,1}(x) = 1
\end{equation}
and if $b \ge 2$ then
\begin{equation}
\label{eqn:h0b}
h_{0,b}(x) = 1 + \frac{x h_{0,b-1}(x)}{1 - x -\ldots - x^{b-1}}.
\end{equation}
\end{proposition}
\begin{proof}
To obtain (\ref{eqn:h01}), observe that $\mu_{0,1} = 1$;  only the empty permutation avoids 1.

To obtain (\ref{eqn:h0b}), first observe that the empty permutation avoids 132, 3241, and $\mu_{0,b}$.
For $n \ge 1$, observe that every permutation $\pi \in S_n(132, 3241, \mu_{0,b})$ for which $\pi(n) = n$ has the form $\tilde{\pi}, n$, where $\tilde{\pi} \in S_{n-1}(132, 3241, \mu_{0,b-1})$, and the resulting map is a bijection between $S_{n-1}(132, 3241, \mu_{0,b-1})$ and the set of permutations in $S_n(132, 3241, \mu_{0,b})$ which end with $n$.
Similarly, observe that if $1 \le i \le b-1$ and $i < n$ then every permutation $\pi \in S_n(132, 3241, \mu_{0,b})$ with $\pi(i) = n$ has the form $n-i+1, n-i+2, \ldots, n, \tilde{\pi}$, where $\tilde{\pi} \in S_{n-i}(132, 3241, \mu_{0,b})$, and the resulting map is a bijection between $S_{n-i}(132, 3241, \mu_{0,b})$ and the set of permutations in $S_n(132, 3241, \mu_{0,b})$ in which $n$ is in the $i$th position.
Finally, observe that if $\pi \in S_n$ avoids 132 and 3241 and $\pi(i) = n$ for $i \ge b$ then $\pi$ does not avoid $\mu_{0,b}$.
Combine these observations to find
$$h_{0,b}(x) = 1 + x h_{0,b-1}(x) + \left( \sum_{i=1}^{b-1} x^i \right) \left( h_{0,b}(x) - 1\right).$$
Solve this equation for $h_{0,b}(x)$ to obtain (\ref{eqn:h0b}).
\end{proof}

We now give a recurrence relation for $h_{a,b}(x)$.

\begin{proposition}
Let $a$ and $b$ denote positive integers.
Then
\begin{equation}
\label{eqn:hab}
h_{a,b}(x) = \frac{1-2x+ xh_{a-1,b}(x)}{(1-x)^2}.
\end{equation}
\end{proposition}
\begin{proof}
First observe that the empty permutation avoids 132, 3241, and $\mu_{a,b}$.
For $n \ge 1$, observe that every permutation $\pi \in S_n(132, 3241, \mu_{a,b})$ for which $\pi(n) = n$ has the form $\tilde{\pi}, n$, where $\tilde{\pi} \in S_{n-1}(132, 3241, \mu_{a,b})$, and the resulting map is a bijection between $S_{n-1}(132, 3241, \mu_{a,b})$ and the set of permutations in $S_n(132, 3241, \mu_{a,b})$ which end with $n$.
Similarly, observe that if $1 \le i \le n-1$ then every permutation $\pi \in S_n(132, 3241, \mu_{a,b})$ with $\pi(i) = n$ has the form $n-i+1, n-i+2, \ldots, n, \tilde{\pi}$, where $\tilde{\pi} \in S_{n-i}(132, 3241, \mu_{a,b})$, and the resulting map is a bijection between $S_{n-i}(132, 3241, \mu_{a-1,b})$ and the set of permutations in $S_n(132, 3241, \mu_{a,b})$ in which $n$ is in the $i$th position.
Combine these observations to find
$$h_{a,b}(x) = 1 + x h_{a,b}(x) + \frac{x}{1-x} (h_{a-1,b}(x) - 1).$$
Solve this equation for $h_{a,b}(x)$ to obtain (\ref{eqn:hab}).
\end{proof}

We now use (\ref{eqn:h0b}) and (\ref{eqn:hab}) to obtain explicit enumerations of $S_n(132, 3241, \mu_{a,b})$ for various small values of $a$ and $b$.

\begin{proposition}
\label{prop:h03}
For all $n \ge 1$,
\begin{equation}
\label{eqn:h03}
|S_n(132, 3241, \mu_{0,3})| = F_{n+2} - 1.
\end{equation}
Moreover,
\begin{equation}
\label{eqn:h03gf}
\sum_{n=0}^\infty |S_n(132, 3241, \mu_{a,b})| x^n = 1 + \frac{x}{(1-x)(1-x-x^2)}.
\end{equation}
\end{proposition}
\begin{proof}
This result is proved in \cite{Mansour33k}, but for completeness we give a proof here.
To obtain (\ref{eqn:h03gf}), first use (\ref{eqn:h01}) and (\ref{eqn:h0b}) to find
$$h_{0,2}(x) = \frac{1}{1-x}.$$
Now set $b = 3$ in (\ref{eqn:h0b}) and use the last equation to simplify the result, obtaining (\ref{eqn:h03gf}).
To obtain (\ref{eqn:h03}), observe that 
$$\frac{x}{(1-x)(1-x-x^2)} = \frac{F_1 x+F_2}{1 - x - x^2} - \frac{1}{1-x}.$$
Now (\ref{eqn:h03}) is immediate from (\ref{eqn:fnpluskgf}) and (\ref{eqn:nplus1choosekgf}).
\end{proof}

We also give a bijective proof of (\ref{eqn:h03}).

\begin{theorem}
\label{thm:h03bijection}
For all $n \ge 1$ there exists a constructive bijection between $S_n(132, 3241, \mu_{0,3})$ and the set of tilings of a $1 \times (n+1)$ rectangle with tiles of size $1 \times 1$ and $1 \times 2$ using at least one $1 \times 2$ tile.
\end{theorem}
\begin{proof}
Suppose we are given such a tiling.
To construct the corresponding permutation, we consider two cases:  the rightmost tile is $1 \times 2$ or the rightmost tile is $1 \times 1$.

If the rightmost tile is $1 \times 2$, then first place a 1 in the rightmost tile.
Now fill the rightmost empty tile with the smallest available numbers, placing them in the tile from left to right in increasing order.
Repeat this process until every tile is filled.

If the rightmost tile is $1 \times 1$, then let $m$ denote the number of $1 \times 1$ tiles to the right of the rightmost $1 \times 2$ tile.
Place an $m$ in the rightmost $1 \times 2$ tile.
Fill the (necessarily $1 \times 1$) tiles to the right of the rightmost $1 \times 2$ with the numbers $m-1, m-2, \ldots, 1, m+1$, in this order.
Now fill the rightmost empty tile with the smallest available numbers, placing them in the tile from left to right in increasing order.
Repeat this process until every tile is filled.

It is routine to verify that this map is invertible, and that the permutation constructed avoids 132, 3241, and $\mu_{0,3}$.
\end{proof}

\begin{example}
Under the bijection given in the proof of Theorem \ref{thm:h03bijection}, the permutation 7563214 corresponds to the tiling $\Box\ \Box\hspace{-2.5pt}\Box\ \Box\hspace{-2.5pt}\Box\ \Box\ \Box\ \Box$ and the permutation 897564321 corresponds to the tiling $\Box\hspace{-2.5pt}\Box\ \Box\ \Box\hspace{-2.5pt}\Box\ \Box\ \Box\ \Box\ \Box\hspace{-2.5pt}\Box$.
\end{example}

We omit the proofs of the remaining propositions in this section, since they are similar to the proof of Proposition \ref{prop:h03}.

\begin{proposition}
For all $n \ge 1$, 
\begin{equation}
|S_n(132, 3241, \mu_{0,4})| = \frac{1}{2} \left( 2 F_{3,n+3} + F_{3,n+2} + F_{3,n} - 2 F_{n+4} + 1\right).
\end{equation}
Moreover,
\begin{displaymath}
\sum_{n=0}^\infty |S_n(132, 3241, \mu_{0,4})| = 1 + \frac{x - x^2 + x^4}{(1-x)(1-x-x^2)(1-x-x^2-x^3)}.
\end{displaymath}
\end{proposition}

\begin{proposition}
For all $n \ge 1$,
\begin{eqnarray*}
\lefteqn{|S_n(132, 3241, \mu_{0,5})| = } \\
 & & \frac{1}{6}\left( 14 F_{4,n+4} + 8 F_{4,n+3} + 4 F_{4,n+2} + 2 F_{4,n} - 6 F_{3,n+6} - 3 F_{3,n+5} - 3 F_{3,n+3} + 6 F_{n+5} - 1 \right).
\end{eqnarray*}
Moreover,
\begin{eqnarray*}
\lefteqn{\sum_{n=0}^\infty |S_n(132, 3241, \mu_{0,5})| x^n = } \\
& & 1 + \frac{x - 2x^2 + 2x^4 + 2x^5 - x^6 - x^7}{(1-x)(1-x-x^2)(1-x-x^2-x^3)(1-x-x^2-x^3-x^4)}.
\end{eqnarray*}
\end{proposition}

\begin{proposition}
For all $n \ge 1$,
\begin{equation}
\label{eqn:omegainverse}
|S_n(132, 3241, \mu_{1,3})| = F_{n+5} - {{n+1} \choose {2}} - 2 {{n+1} \choose {1}} - 2.
\end{equation}
Moreover,
\begin{displaymath}
\sum_{n=0}^\infty |S_n(132, 3241, \mu_{1,3})| x^n = \frac{1 - 3x + 3x^2 + x^3 - x^4}{(1-x)^3 (1 - x -x^2)}.
\end{displaymath}
\end{proposition}

Comparing (\ref{eqn:omegainverse}) with (\ref{eqn:omega4}) we find that for $n \ge 1$,
$$|S_n(132, 3241, \mu_{1,3})| = |S_n(132, 2341, \omega_4)|.$$
This is no surprise, however, since $S_n(132, 2341, \omega_4)$ is precisely the set of inverses of the elements of $S_n(132, 3241, \mu_{1,3})$.

\begin{proposition}
For all $n \ge 1$,
\begin{equation}
|S_n(132, 3241, \mu_{2,3})| = F_{n+8} - {{n+1}\choose{4}} - 3 {{n+1} \choose {3}} - 4 {{n+1} \choose {2}} - 9 {{n+1} \choose {1}} - 11.
\end{equation}
Moreover,
\begin{displaymath}
\sum_{n=0}^\infty |S_n(132, 3241, \mu_{2,3})| x^n = \frac{1 - 5x+10x^2-8x^3+x^4+4x^5-2x^6}{(1-x)^5 (1-x-x^2)}.
\end{displaymath}
\end{proposition}

\begin{proposition}
For all $n \ge 1$,
\begin{equation}
|S_n(132, 3241, \mu_{1,4})| = \frac{1}{2}\left( F_{3,n+5} + 2 F_{3,n+4} \right) - F_{n+7} + \frac{1}{2} {{n+1}\choose{2}} + \frac{3}{2} {{n+1} \choose {1}} + 5.
\end{equation}
Moreover,
\begin{displaymath}
\sum_{n=0}^\infty |S_n(132, 3241, \mu_{1,4})| x^n = \frac{1-4x+5x^2-x^4-x^5+x^7}{(1-x)^3(1-x-x^2)(1-x-x^2-x^3)}.
\end{displaymath}
\end{proposition}

\section{Extended Restrictions}

Let $R$ denote a set of permutations.
In this section we describe a method of building sets of forbidden subsequences from $R$ for which the associated sets of restricted permutations can be easily enumerated in terms of $|S_n(R)|$.
We then use this method to produce some specific enumerations, several of which involve Fibonacci numbers or $k$-generalized Fibonacci numbers.
We begin by describing our technique for building on $R$.

\begin{definition}
Suppose $\alpha \in S_{n-1}$ and $\beta \in S_n$.
We say $\beta$ is an {\upshape extension} of $\alpha$ whenever $\alpha$ is the permutation obtained by removing $n$ from $\beta$.
We observe that every permutation in $S_n$ has exactly $n+1$ extensions.
\end{definition}

\begin{definition}
Let $R$ denote a set of permutations.
We write $E(R)$ to denote the set of permutations $\pi$ such that $\pi$ is an extension of an element of $R$, and we refer to $E(R)$ as the {\upshape extension of $R$}.
More generally, we write $E^0(R) = R$ and for any $k \ge 2$ we define $E^k(R)$ inductively, so that $E^k(R) = E(E^{k-1}(R))$.
\end{definition}

Our results in this section are based on the following observation.

\begin{theorem}
Let $R$ denote a set of permutations.
Then for all $n \ge 1$,
\begin{equation}
\label{eqn:SEES}
S_n(E(R)) = E(S_{n-1}(R)).
\end{equation}
Moreover, for all $k \ge 0$,
\begin{equation}
\label{eqn:nkSn}
|S_n(E^k(R))| = \frac{n!}{(n-k)!} |S_{n-k}(R)| \hspace{30pt} (n \ge k)
\end{equation}
and
\begin{equation}
\label{eqn:n!}
|S_n(E^k(R))| = n! \hspace{30pt} (0 \le n < k).
\end{equation}
\end{theorem}
\begin{proof}
To obtain (\ref{eqn:SEES}), suppose $\pi' \in S_{n-1}$, $\pi$ is an extension of $\pi'$, and $\sigma' \in S_m$.
Observe that $\pi'$ contains a subsequence of type $\sigma'$ if and only if there exists an extension $\sigma$ of $\sigma'$ such that $\pi$ contains a subsequence of type $\sigma$.
Now (\ref{eqn:SEES}) is immediate.

To obtain (\ref{eqn:nkSn}), we argue by induction on $k$.
Since $E^0(R) = R$, the result is immediate for $k = 0$.
Now observe that every permutation in $S_n$ has exactly $n+1$ extensions, and every permutation in $S_{n+1}$ is the extension of exactly one permutation in $S_n$.
Therefore, using (\ref{eqn:SEES}) and induction we find
\begin{eqnarray*}
|S_n(E^k(R))| &=& |E(S_{n-1}(E^{k-1}(R)))| \\
&=& n |S_{n-1}(E^{k-1}(R))| \\
&=& \frac{n!}{(n-k)!} |S_{n-k}(R)|,
\end{eqnarray*}
as desired.

To obtain (\ref{eqn:n!}), observe that every permutation in $R$ has length at least $k$, so $S_n(R) = S_n$ when $n < k$.
\end{proof}

We now give several applications of this theorem.

\begin{proposition} 
For all $k \ge 0$ and all $n \ge k$,
$$|S_n(E^k(123,132,213))| = \frac{n!}{(n-k)!} F_{n+1-k}.$$
In particular, for all $n \ge 1$,
$$|S_n(1234,1243,1423,4123,1324,1342,1432,4132,2134,2143,2413,4213)| = n F_n.$$	
\end{proposition}
\begin{proof}
Set $R = \{123, 132, 213\}$ in (\ref{eqn:nkSn}) and use (\ref{eqn:SimionSchmidtintro}) to simplify the result.
\end{proof}

\begin{proposition}
Let $R$ denote one of the following sets of permutations.
$$\begin{array}{c}
\{123,1432\};\ \{123,2143\};\ \{123,2413\};\ \{132,1234\};\ \{132,2134\};\\
\{132,2314\};\ \{132,2341\};\ \{132,3241\};\ \{132,3412\}
\end{array}$$
For all $k \ge 0$ and all $n \ge k$,
\begin{equation}
\label{eqn:nWest}
|S_n(E^k(R))| = \frac{n!}{(n-k)!} F_{2(n-k)-1}.
\end{equation} 
\end{proposition}
\begin{proof}
In \cite{WestGenTrees} West shows that for each of the sets listed, $S_n(R) = F_{2n-1}$.
Combine this with (\ref{eqn:nkSn}) to obtain (\ref{eqn:nWest}).
\end{proof}

\begin{proposition} 
(see Guibert \cite{G})
Fix $\tau \in S_3$.
For all $k \ge 0$ and all $n \ge k$,
\begin{equation}
\label{eqn:nCn}
|S_n(E^k(\tau))| = \frac{n!}{(n-k+1)!}{{2n-2k}\choose{n-k}}.
\end{equation} 
\end{proposition}
\begin{proof}
It is well known that for all $\tau \in S_3$ we have $S_n(\tau) = C_n$, where ${\displaystyle C_n = \frac{1}{n+1} {{2n} \choose {n}}}$ is the $n$th Catalan number.
Combine this with (\ref{eqn:nkSn}) to obtain (\ref{eqn:nCn}).
\end{proof}

\section{Generalizations of Simion and Schmidt Involving Larger Sets of Restrictions}

Let $R$ denote a set of permutations.
Up to this point we have singled out those $R$ for which $|S_n(R)|$ can be expressed in terms of Fibonacci numbers or $k$-generalized  Fibonacci numbers.
Put another way, we have considered those $R$ for which $|S_n(R)|$ is closely related to certain linear homogeneous recurrence relations with constant coefficients.
In this section we describe a family of sets of permutations such that if $R$ is a set in the family then $|S_n(R)|$ satisfies a linear homogeneous recurrence relation with constant coefficients.
We begin with some terminology.

\begin{definition}
Let $k$ and $l$ denote positive integers such that $l \le k$ and let $a_1, \ldots, a_l$ denote a sequence of positive integers.
We say a permutation $\sigma \in S_k$ {\upshape agrees with $a_1, \ldots, a_l$ to length $j$} whenever $\sigma(i) = a_i$ for $1 \le i \le j$.
If $\sigma(1) \neq a_1$ then we say $\sigma$ agrees with $a_1, \ldots, a_l$ to length 0.
\end{definition}

We now define our sets of forbidden subsequences.

\begin{definition}
Let $k$ and $l$ denote positive integers such that $l \le k$ and let $a_1, \ldots, a_l$ denote a sequence of positive integers.
We write $R^k_{a_1, \ldots, a_l}$ to denote the set of permutations in $S_k$ such that $\sigma \in R^k_{a_1, \ldots, a_l}$ whenever one of the following holds.
\begin{enumerate}
\item $\sigma$ agrees with $a_1, \ldots, a_l$ to length $l$.
\item There exists $i$, $0 \le i \le l-1$, such that $\sigma$ agrees with $a_1, \ldots, a_l$ to length $i$ and $\sigma(i+1) < a_{i+1}$.
\end{enumerate}
\end{definition}

The sequence $|S_n(R^k_{a_1, \ldots, a_l})|$ satisfies a certain linear homogeneous recurrence relation with constant coefficients.

\begin{theorem}
\label{thm:tt}
Let $k$ and $l$ denote positive integers such that $l \le k$ and let $a_1, \ldots, a_l$ denote a sequence of positive integers in which $a_1, \ldots, a_{l-1}$ are distinct.
Then for all $n < k$ we have
\begin{equation}
\label{eqn:kln!}
|S_n(R^k_{a_1, \ldots, a_l})| = n!
\end{equation}
and for all $n \ge k$ we have
\begin{equation}
\label{eqn:klrecurrence}
|S_n(R^k_{a_1, \ldots, a_l})| = \sum_{j=1}^l (k-a_j-\eta_j) |S_{n-j}(R^k_{a_1, \ldots, a_l})|,
\end{equation}
where $\eta_j = \{i\ |\ a_i > a_j\ \mbox{and}\ i < j\}$ for $1 \le j \le l$.
\end{theorem}
\begin{proof}
To obtain (\ref{eqn:kln!}), observe that $R^k_{a_1, \ldots, a_l} \subseteq S_k$, so $S_n(R^k_{a_1, \ldots, a_l}) = S_n$ when $n < k$.

To prove (\ref{eqn:klrecurrence}), we first set some notation.
For any $j$, $0 \le j \le l$, let $P_{n,j}$ denote the set of permutations $\pi \in R^k_{a_1, \ldots, a_l}$ such that $\pi(i) = n - k + a_i$ for $1 \le i \le j$.
We observe that $P_{n,0} = S_n(R^k_{a_1, \ldots, a_l})$.
We also observe that if $\pi \in S_n$ satisfies $\pi(i) = n - k + a_i$ for $ 1 \le i \le l$ then the subsequence $n-k+a_1, n-k+a_2, \ldots, n-k+a_l$ of $\pi$, together with the remaining elements of $\pi$ which are greater than $n-k$, forms a pattern in $R^k_{a_1, \ldots, a_l}$.
Therefore $P_{n,l}$ is empty.

We make the following claim concerning $P_{n,j}$.

\medskip
\noindent
{\bf Claim One}
For all $j$, $0 \le j \le l-1$,
\begin{equation}
\label{eqn:claimone}
|P_{n,j}| =  (k - a_{j+1} - \eta_{j+1}) |S_{n-j-1}(R)| + |P_{n,j+1}|,
\end{equation}
where $\eta_j = \{i\ |\ a_i > a_j\ \mbox{and}\ i < j\}$ for $1 \le j \le l$.

\medskip
Assuming for the moment that claim one holds, we complete the proof of (\ref{eqn:klrecurrence}).
Set $j = 0$ in (\ref{eqn:claimone}) and recall that $P_{n,0} = S_n(R^k_{a_1, \ldots, a_l})$ to find
$$|S_n(R^k_{a_1, \ldots, a_l})| = (k - a_1 - \eta_1) |S_{n-1}(R^k_{a_1, \ldots, a_l})| + |P_{n,1}|.$$
Now use (\ref{eqn:claimone}) with $j = 1$ to eliminate $|P_{n,1}|$ in the last equation.
Repeat the process to obtain
$$|S_n(R^k_{a_1, \ldots, a_l})| = \sum_{j=1}^l (k - a_j - \eta_j) |S_{n-j}(R^k_{a_1, \ldots, a_l})| + |P_{n,l}|.$$
Since $P_{n,l}$ is empty, (\ref{eqn:klrecurrence}) follows.

\medskip

\noindent
{\it Proof of Claim One.}
First observe that $P_{n,j+1} \subseteq P_{n,j}$, so we restrict our attention to those $\pi \in P_{n,j}$ for which $\pi(j+1) \neq n - k + a_{j+1}$.
We consider two cases:  $\pi(j+1) < n - k + a_{j+1}$ and $\pi(j+1) > n - k + a_{j+1}$.

Suppose $\pi(j+1) < n - k + a_{j+1}$.
Then the subsequence $n - k + a_1, n - k + a_2, \ldots, n - k + a_j, \pi(j+1)$ of $\pi$, together with the remaining elements of $\pi$ which are greater than $n - k + 1$, forms a pattern in $R_{a_1, \ldots, a_l}^k$.
This contradicts the fact that $\pi$ avoids every pattern in $R_{a_1, \ldots, a_l}^k$, so there is no $\pi \in P_{n,j}$ for which $\pi(j+1) < n - k + a_{j+1}$.

Suppose $\pi(j+1) > n - k + a_{j+1}$.
Observe that if there exists $i$, $1 \le i \le j$, for which $\pi(i)$ participates in a pattern in $R^k_{a_1, \ldots, a_k}$, then that pattern must involve all of those entries of $\pi$ which are greater than $n-k$.
For instance, that pattern must involve $\pi(1), \pi(2), \ldots, \pi(j)$, since each of these is greater than $n-k$.
Similarly, such a pattern must involve $\pi(j+1)$, since $\pi(j+1) > n-k+a_{j+1} > n-k$.
However, because $\pi(j+1) > n-k+a_{j+1}$, there is no pattern in $R^k_{a_1, \ldots, a_l}$ whose first $j+1$ entries have the same relative order as $\pi(1), \ldots, \pi(j+1)$.
Therefore if $\pi \in P_{n,j}$ and $\pi(j+1) > n-k+a_{j+1}$ then $\pi$ has the form $n-k+a_1, n-k+a_2, \ldots, n-k+a_j, r, \tilde{\pi}$, where $n-k+a_{j+1} < r \le n$ and $\tilde{\pi}$ is a permutation of the remaining numbers which avoids $R^k_{a_1, \ldots, a_l}$.
Since there are $k-a_{j+1} - \eta_{j+1}$ choices for $r$, and $|S_{n-j-1}(R_{a_1, \ldots, a_l})|$ choices for $\tilde{\pi}$, the number of permutations $\pi \in P_{n,j}$ for which $\pi(j+1) > n-k+a_{j+1}$ is equal to $(k-a_{j+1}- \eta_{j+1}) |S_{n-j-1}(R_{a_1, \ldots, a_l})|$.
Now (\ref{eqn:claimone}) follows.
\end{proof}

Simion and Schmidt's result (\ref{eqn:SimionSchmidtintro}) is a special case of Theorem \ref{thm:tt}.
To see this, first observe that $R^3_{2,2} = \{123, 132, 213\}$.
In view of this, we set $k = 3$, $l = 2$, $a_1 = 2$, and $a_2 = 2$ in (\ref{eqn:klrecurrence}) to obtain
$$|S_n(123, 132, 213)| = |S_{n-1}(123, 132, 213)| + |S_{n-2}(123, 132, 213)|$$
for $n \ge 3$.
Since $F_2 = 1!$ and $F_3 = 2!$, we find that $|S_n(123, 132, 213)| = F_{n+1}$ for $n \ge 1$.

Theorem \ref{thm:tt} has many other interesting special cases.
We conclude this section by mentioning just a few of them.

\begin{proposition}
(\cite[Thm. 1]{Mproc})
Let $k$ denote a positive integer and fix $a$ such that $1 \le a \le k$.
Then for all $n \ge k-1$, 
\begin{equation}
\label{eqn:powerfact}
|S_n(R^k_a)| = (k-1)! (k-a)^{n-k+1}.
\end{equation}
\end{proposition}
\begin{proof}
Set $l = 1$ and $a_1 = a$ in (\ref{eqn:klrecurrence}) to find
$$|S_n(R^k_a)| = (k-a) |S_{n-1}(R^k_a)|$$
for $n \ge k$.
Combine this with (\ref{eqn:kln!}) to obtain (\ref{eqn:powerfact}).
\end{proof}

\begin{proposition}
For all $n \ge 1$,
$$|S_n(1234,1243,1324,1342,1423,1432,2134,2143,2314,2341)| = p_n + p_{n-2},$$
where $p_n$ is the $n$th Pell number, which is defined by $p_0 = 0$, $p_1 = 1$, and $p_n = 2 p_{n-1} + p_{n-2}$ for $n \ge 2$.
\end{proposition}
\begin{proof}
For notational convenience, we abbreviate 
$$b_n = |S_n(1234,1243,1324,1342,1423,1432,2134,2143,2314,2341)|.$$
Setting $l = 2$, $k = 4$, $a_1 = 2$, and $a_2 = 3$ in Theorem \ref{thm:tt} we find that $b_n = 2 b_{n-1} + b_{n-2}$ for $n \ge 4$.
Since $p_n + p_{n-2} = n!$ for $n = 1, 2, 3$, the result follows.
\end{proof}

For more information on the Pell numbers, see \cite{FibProblem}, \cite[pp. 122--125]{Beiler}, \cite{Emerson}, and \cite{Horadam}.

We conclude the paper with one more generalization of Simion and Schmidt's result (\ref{eqn:SimionSchmidtintro}).

\begin{proposition}
Let $k$ denote a positive integer such that $k \ge 3$.
Then for all $n \ge k$,
\begin{displaymath}
|S_n(R^k_{k-1,k-1})| = (k-2)! ( F_{n-k+4} + (k-3) F_{n-k+2} ).
\end{displaymath}
\end{proposition}
\begin{proof}
Consider the sequence $x_{k-2}, x_{k-1}, \ldots$ in which
$$x_n = |S_n(R^k_{k-1,k-1})|.$$
By (\ref{eqn:kln!}) the first two terms in this sequence are $(k-2)!$ and $(k-1)!$.
By (\ref{eqn:klrecurrence}) this sequence satisfies $x_n = x_{n-1} + x_{n-2}$ for all $n \ge k$.
Now consider the sequence $y_{k-2}, y_{k-1}, \ldots$ in which $y_{k-2} = (k-2)!$, $y_{k-1} = (k-1)!$, and $y_n = (k-2)! ( F_{n-k+4} + (k-3) F_{n-k+2} )$ for $n \ge k$.
Clearly $x_{k-2} = y_{k-2}$ and $x_{k-1} = y_{k-1}$.
Moreover, it is routine to verify that $y_n = y_{n-1} + y_{n-2}$ for all $n \ge k$.
Therefore $y_n = x_n$ for all $n \ge k$, as desired.
\end{proof}


\end{document}